\documentclass[12pt]{article}
\usepackage{booktabs}
\usepackage{hyperref}
\hypersetup{
 colorlinks=true,
 anchorcolor=yellow,
 linkcolor=cyan,
 filecolor=blue,
 urlcolor=red,
 citecolor=green,
}
\usepackage{caption}
\usepackage{mathrsfs}
\usepackage{amsmath}
\usepackage{amsfonts,amsthm,amssymb,mathrsfs,bbding}
\usepackage{float}

\usepackage{graphicx}
\usepackage{enumerate}
\usepackage{tikz}
\usepackage{caption}
\usepackage{rotating}
\allowdisplaybreaks[4]
\usepackage{tabularx}
\usepackage{cite}
\evensidemargin=\oddsidemargin\topmargin=-1.5cm

\newtheorem{thm}{Theorem}[section]
\newtheorem{prob}{Problem}[section]

\newtheorem{lem}{Lemma}[section]
\newtheorem{cor}{Corollary}[section]

\newtheorem{conj}{Conjecture}[section]

\theoremstyle{definition}
\renewcommand\proofname{\bf Proof}
\addtocounter{section}{0}
\begin{document}

\title{\LARGE{\bf Spectral radius and $[a,b]$-factors in graphs} \footnote{This work is supported by the National Natural Science Foundation of China (Grant Nos. 11771141, 12011530064 and 11871391)}\setcounter{footnote}{-1}\footnote{\emph{Email address:} ddfan0526@163.com (D. Fan), huiqiulin@126.com (H. Lin), luhongliang@mail.xjtu.edu.cn(H.Lu).}}
\author{Dandan Fan$^a$,  Huiqiu Lin$^a$\thanks{Corresponding author.}~~and Hongliang Lu$^b$\\[2mm]
\small\it $^a$ School of Mathematics, East China University of Science and Technology, \\
\small\it   Shanghai 200237, P.R. China\\[1mm]
\small\it $^b$ School of Mathematics and Statistics, Xi'an Jiaotong University\\
\small\it Xi'an, Shanxi 710049, China \\}
\date{}
\maketitle
{\flushleft\large\bf Abstract}
An $[a,b]$-factor of a graph $G$ is a spanning subgraph $H$ such that $a\leq d_{H}(v)\leq b$ for each $v\in V(G)$. In this paper, we provide spectral conditions for the existence of an odd $[1,b]$-factor in a connected graph with minimum degree $\delta$ and the existence of an $[a,b]$-factor in a graph, respectively. Our results generalize and improve some previous results on perfect matchings of graphs. For $a=1$, we extend the result of O\cite{S.O} to obtain an odd $[1,b]$-factor and further improve the result of Liu, Liu and Feng\cite{W.L} for $a=b=1$. For $n\geq 3a+b-1$, we confirm the conjecture of Cho, Hyun, O and Park\cite{E.C}. We conclude some open problems in the end.

\begin{flushleft}
\end{flushleft}

\section{Introduction}

The study of factors from graph eigenvalues has a rich history. In 2005, Brouwer and Haemers \cite{A.B} described a regular graph to contain a perfect matching in terms of the third largest eigenvalue. Their result was improved in \cite{S.C-1,S.C,S.C-2} and extended in \cite{H.L-3,H.L-4} to obtain a regular factor. Very recently, O \cite{S.O} provided a spectral condition to guarantee the existence of a perfect matching in a graph. One problem concerning perfect matching which has attracted considerable interest is that of determining the structure of graphs with a unique perfect matching\cite{X.W,D.B,S.P}. Yang and Ye\cite{Y.Y} characterized all bipartite graphs with a unique perfect matching whose adjacency matrices have inverses diagonally similar to non-negative matrices, which settles an open problem of Godsil on inverses of bipartite graphs in \cite{C.G}. Lov\'{a}sz\cite{L.L} proved that a graph of order $n$ with a unique perfect matching cannot have more than $\frac{n^2}{4}$ edges. Based on the structural properties of graphs with a unique perfect matching, in this paper, we first determine the graph attaining the maximum spectral radius among all graphs of order $2n$ with a unique perfect matching.

Let $G$ be a graph with adjacency matrix $A(G)$. The largest eigenvalue of $A(G)$, denoted by $\rho(G)$, is called the \textit{spectral radius} of $G$. Denote by $'\nabla'$ and $'\cup'$ the join and union products, respectively. Suppose that $G_1$ is an empty graph with vertex set $U=\{u_{1},u_{2},\ldots, u_{n}\}$, and $G_2$ is a complete graph with vertex set $W=\{w_{1},w_{2},\ldots, w_{n}\}$. Let $G(2n,1)$ be the graph of order $2n$ obtained from $G_{1}\cup G_{2}$ by letting $N_{G_2}(u_i)=\{w_{1},w_{2},\ldots,w_{i}\}$ for $1\leq i\leq n$. Clearly, $G(2n,1)$ contains a unique perfect matching.
\begin{thm}\label{thm::1.1}
If $G$ is a connected graph of order $2n$ with a unique perfect matching, then $\rho(G)\leq\rho(G(2n,1))$, with equality if and only if $G\cong G(2n,1)$.
\end{thm}

Another problem concerning perfect matching is extending it to a general factor. An $\textit{odd [1,b]-factor}$ of a graph $G$ is a spanning subgraph $H$ such that $d_{H}(v)$ is odd and $1\leq d_{H}(v)\leq b$ for each $v\in V(G)$. Lu, Wu and Yang \cite{H.L} gave a sufficient condition for the existence of an odd $[1,b]$-factor in a graph in terms of the third largest eigenvalue. Recently, Kim, O, Park and Ree \cite{S.K} improved the spectral condition of Lu, Wu and Yang. In this paper, we generalize their result by giving a spectral condition to guarantee the existence of an odd $[1,b]$-factor with minimum degree $\delta$. Let $T(n,b,\delta)=K_{\delta} \nabla (K_{n-(b+1)\delta-1}\cup (b\delta+1)K_{1})$ and $F(b,\delta)=\max\{4(b+1)\delta+4, b\delta^3+\delta\}$, where $b$ is a positive odd integer.
\begin{thm}\label{thm::1.2}
Suppose that $G$ is a connected graph of even order $n\geq F(b,\delta)$ with minimum degree $\delta$. If $\rho(G)\geq\rho(T(n,b,\delta))$, then $G$ contains an odd $[1,b]$-factor, unless $G\cong T(n,b,\delta)$.
\end{thm}

An \textit{$[a,b]$-factor} of a graph $G$ is a spanning subgraph $H$ such that $a\leq d_{H}(v)\leq b$ for each $v\in V(G)$. In the past decades, some researchers provided various parameter conditions for a graph to have an $[a,b]$-factor, such as the degree condition\cite{Y.L}, the neighborhood condition \cite{J.L-1}, the stability number\cite{M.K}, and the binding number\cite{C.C}. In addition, the conditions for a graph to have a fractional $[a,b]$-factor (see \cite{E.S} for the definition) was also investigated by several researchers\cite{Y.E,H.L-2,M.K-1}. Only very recently, Cho, Hyun, O and Park\cite{E.C} posed a conjecture regarding the spectral condition for the existence of an $[a,b]$-factor as follows.

\begin{conj}(See \cite{E.C})\label{conj::1.1}
Let $a\cdot n$ be an even integer at least 2, where $n\geq a+1$. If $G$ is a graph of order $n$ with $\rho(G)>\rho(H_{n,a})$ where $H_{n,a}=K_{a-1}\nabla(K_{1}\cup K_{n-a})$, then $G$ contains an $[a,b]$-factor.
\end{conj}

In this paper, we confirm Conjecture \ref{conj::1.1}  for $n\geq 3a+b-1$.
\begin{thm}\label{thm::1.3}
 Let $a\cdot n$ be an even integer, where $n\geq 3a+b-1$ and $b\geq a\geq 1$. If $\rho(G)>\rho(H_{n,a})$, then $G$ contains an $[a,b]$-factor.
\end{thm}

A $[k,k]$-factor is called a $k$-factor. As a corollary of Theorem \ref{thm::1.3}, we have the following result.
\begin{cor}\label{cor::1.1}
Let $k\cdot n$ be an even integer, where $n\geq 4k-1$. If $\rho(G)>\rho(H_{n,k})$, then $G$ contains a $k$-factor.
\end{cor}

\section{Proof of Theorem \ref{thm::1.1}}
If $G$ is connected, then $A(G)$ is irreducible. By the Perron-Frobenius theorem(cf. \cite[Section 8.8]{C.G-2}), the Perron vector $x$ is a positive eigenvector of $A(G)$ respect to $\rho(G)$. For any $v\in V(G)$, let $N_{G}(v)$ and $d_G(v)$ be the neighborhood and degree of $v$, respectively. In this section, we give the proof of Theorem \ref{thm::1.1}. Before proceeding, the following lemmas are needed.
\begin{lem}(See \cite{H.L-1})\label{lem::2.1}
Let $G$ be a connected graph, and let $u,v$ be two vertices of $G$. Suppose that $v_{1},v_{2},\ldots,v_{s}\in N_{G}(v)\backslash N_{G}(u)$ with $s\geq 1$, and $G^*$ is the graph obtained from $G$ by deleting the edges $vv_{i}$ and adding the edges $uv_{i}$ for $1\leq i\leq s$. Let $x$ be the Perron vector of $A(G)$. If $x_{u}\geq x_{v}$, then $\rho(G)<\rho(G^*)$.
\end{lem}

An edge is said to be a \emph{cut edge} if its removal increases the number of components of a graph.
\begin{lem}(See \cite{A.K})\label{lem::2.2}
Let $G$ be a connected graph with a unique perfect matching. Then $G$ contains a cut edge $uv$ that is an edge of the perfect matching of $G$.
\end{lem}

\begin{lem}(See \cite{BH})\label{lem::2.3}
If $u$ and $v$ are two nonadjacent vertices of graph $G$, then $\rho(G+uv)>\rho(G)$.
\end{lem}

For any $S\subseteq V(G)$, let $G[S]$ be the subgraph of $G$ induced by $S$. An edge of a graph is called \emph{pendant edge} if exactly one of its ends has the degree 1. We now give a short proof of Theorem \ref{thm::1.1}.

\renewcommand\proofname{\bf Proof of Theorem \ref{thm::1.1}}
\begin{proof}

Suppose that $G$ is a connected graph of order $2n$ with a unique perfect matching $M$. By Lemma \ref{lem::2.2}, there exists a cut edge, say $u_{0}v_{0}$, that is contained in $M$. Observe   that $G-u_{0}v_{0}$ consists of two odd components. Let $x^{(0)}$ be the Perron vector of $A(G)$. Without loss of generality, we assume that $x^{(0)}_{u_0}\geq x^{(0)}_{v_0}$. Let
\begin{equation*}
G_1= G-\sum_{w\in N_{G}(v_0)\setminus\{u_0\}}v_{0}w+\sum_{w\in V(G)\setminus(N(u_0)\cup \{u_0\})}u_{0}w.
\end{equation*}
We see that $G_1$ also has a unique perfect matching, say $M_1$. By Lemmas \ref{lem::2.1} and \ref{lem::2.3}, we have $\rho(G)\leq\rho(G_1)$, with equality  if and only if $G\cong G_{1}$. Let $S_{1}=V(G_1)-\{u_{0},v_{0}\}$. Since $u_{0}v_{0}$ is a pendant edge of $G_1$ that is contained in $M_1$, the induced subgraph $G_1[S_1]$ also contains a unique perfect matching, i.e., $M_1\setminus\{u_0v_0\}$. Again by Lemma \ref{lem::2.2}, there exists a cut edge $u_{1}v_{1}$ in $G_1[S_1]$ that is contained in $M_1\setminus\{u_0v_0\}$. Let $x^{(1)}$ be the Perron vector of $A(G_1)$. Assume that $x^{(1)}_{u_1}\geq x^{(1)}_{v_1}$. Let
\begin{equation*}
G_2= G_1-\sum_{w\in N_{G_1[S_1]}(v_1)\setminus\{u_1\}}v_{1}w+\sum_{w\in S_1\setminus(N(u_1)\cup \{u_1\})}u_{1}w.
\end{equation*}
Clearly, $G_2$ also has a unique perfect matching. Again by Lemmas \ref{lem::2.1} and \ref{lem::2.3},  $\rho(G_1)\leq\rho(G_2)$, where the equality holds if and only if $G_{1}\cong G_{2}$. By repeating this procedure, we can construct a  sequence of graphs $G_{0}, G_{1}, G_{2},\cdots, G_{n-1}$:
\begin{itemize}
\item[(i)] $G_0=G$;
\item[(ii)] for $i\in [0,n-2]$, let $S_{i}=V(G_{i})-\{v_0,v_1,\ldots,v_{i-1},u_0,u_1,\ldots,u_{i-1}\}$ and
\begin{equation*}
G_{i+1}= G_{i}-\sum_{w\in N_{G_i[S_i]}(v_i)\setminus\{u_i\}}v_{i}w+\sum_{w\in S_i\setminus(N(u_i)\cup \{u_i\})}u_{i}w,
\end{equation*}
where $u_iv_i$ is a cut edge of $G_i[S_i]$ that is contained in the unique perfect matching of $G_i[S_i]$ and  $x^{(i)}_{v_i}\leq x^{(i)}_{u_i}$, where $x^{(i)}$ is the Perron vector of $A(G_i)$.
\end{itemize}
As above, we see that $G_{i}$ has a unique perfect matching for each $i$, and $\rho(G_{i})\leq \rho(G_{i+1})$ with equality if and only if $G_{i}\cong G_{i+1}$  ($0\leq i\leq n-2$). Note that $G_{n-1}\cong G(2n,1)$. Thus we conclude that $\rho(G)\leq\rho(G(2n,1))$, where the equality holds if and only if $G\cong G(2n,1)$.

We complete the proof.
\end{proof}

\section{Proof of Theorem \ref{thm::1.2}}
The well-known sufficient and necessary condition for the existence of an odd $[1,b]$-factor established by Amahashi\cite{A.A}.
\begin{lem}(See \cite{A.A})\label{lem::3.1}
Let $G$ be a graph and let $b$ be a positive odd integer. Then $G$ contains an odd $[1,b]$-factor if and only if for every subset $S\subseteq V(G)$,
$$o(G-S)\leq b|S|,$$
where $o(H)$ is the number of odd components in a graph $H$.
\end{lem}

\begin{lem}(See \cite{D.F})\label{lem::3.2}
Let $n=\sum_{i=1}^t n_i+s$. If $n_{1}\geq n_{2}\geq \cdots\geq n_{t}\geq p$ and $n_{1}<n-s-p(t-1)$, then
$$\rho(K_{s} \nabla (K_{n_{1}}\cup K_{n_{2}}\cup \cdots \cup K_{n_{t}}))<\rho(K_{s} \nabla (K_{n-s-p(t-1)}\cup (t-1)K_{p})).$$
\end{lem}
Now we shall give the proof of Theorem \ref{thm::1.2}.
\renewcommand\proofname{\bf Proof of Theorem \ref{thm::1.2}}
\begin{proof}

Suppose that $G$ contains no odd $[1,b]$-factor, by Lemma \ref{lem::3.1}, there exists some nonempty subset $S$ of $V(G)$ such that $q=o(G-S)> b|S|$.
Since $n$ is even, $q$ and $b|S|$ have the same parity,
we have $q\geq b|S|+2$. Let $|S|=s$ and $t=bs+2$. Then $G$ is a spanning subgraph of $G_s^1=K_{s} \nabla (K_{n_1}\cup K_{n_2}\cup \cdots \cup K_{n_t})$ for some positive odd integers $n_1\geq n_2\geq \cdots\geq n_t$ with  $\sum_{i=1}^{t}n_i=n-s$. Thus,
\begin{equation}\label{equ::1}
\rho(G)\leq\rho(G_s^1),
\end{equation}
where the equality holds if and only if $G\cong G_s^1$. Note that the graph $G_s^1$ contains no odd $[1,b]$-factor. Then we shall derive the proof into the following three cases.

{\flushleft\bf{Case 1.}} $s\geq\delta+1.$

Let $G_s^2=K_{s} \nabla (K_{n-s-t+1}\cup (t-1)K_{1})$. Note that $s\geq \delta+1$ and  $t=bs+2$. Then by Lemma \ref{lem::3.2}, we have
\begin{equation}\label{equ::2}
\rho(G_s^1)\leq \rho(G_{s}^2),
\end{equation}
with equality if and only if $(n_1,\ldots,n_t)=(n-s-t+1,1,\ldots,1)$.
The vertex set of $G_{s}^2$ can be partitioned as $V(G_{s}^2)=V(K_{s})\cup V((bs+1)K_{1})\cup V(K_{n-(b+1)s-1})$, where $V(K_s)=\{v_1,\ldots,v_s\}$, $V((bs+1)K_{1})=\{u_{1},\ldots,u_{bs+1}\}$ and $V(K_{n-(b+1)s-1})=\{w_{1},\ldots,w_{n-(b+1)s-1}\}$. Let
$$G'_s=G_{s}^2+\sum_{i=b\delta+2}^{bs+1}\sum_{j=1}^{n-(b+1)s-1}u_{i}w_{j}+
\sum_{i=b\delta+2}^{bs}\sum_{j=i+1}^{bs+1}u_{i}u_{j}-
\sum_{i=\delta+1}^{s}\sum_{j=1}^{b\delta+1}v_{i}u_{j}.$$
Clearly, $G'_s\cong K_{\delta} \nabla (K_{n-(b+1)\delta-1}\cup (b\delta+1)K_{1})$.
Let $x$ be the Perron vector of $A(G_{s}^2)$ with respect to $\rho=\rho(G_{s}^2)$. By symmetry, $x$ takes the same value (say $x_1$, $x_2$ and $x_3$) on the vertices of $V(K_s)$, $V((bs+1)K_{1})$ and $V(K_{n-(b+1)s-1)})$, respectively. Then, by $A(G_{s}^2)x=\rho x$, we have
\begin{eqnarray*}
 &\rho x_2&=s x_{1},\\
 &\rho x_3 &=s x_{1}+(n-(b+1)s-2) x_{3}.
\end{eqnarray*}
Note that $n\geq s+t= (b+1)s+2$. Then $x_{3}\geq x_{2}$ and
\begin{eqnarray}
 x_{2}=\frac{sx_1}{\rho}. \label{equ::3}
\end{eqnarray}
Similarly, let $y$ be the Perron vector of $A(G'_s)$ corresponding to $\rho(G'_s)=\rho'$. By symmetry, $y$ takes the same values $y_{1}$, $y_{2}$ and $y_{3}$ on the vertices of $V(K_\delta)$, $V((b\delta+1)K_{1})$ and $V(K_{n-(b+1)\delta-1})$, respectively. Then, by $A(G'_s)y=\rho' y$, we obtain
\begin{eqnarray}
 &\rho' y_2&=\delta y_{1},\label{equ::4}\\
 &\rho' y_3 &=\delta y_{1}+(n-(b+1)\delta-2) y_{3}.\label{equ::5}
\end{eqnarray}
Note that $G'_s$ contains $K_{n-(b+1)\delta-1}\cup K_{\delta}\cup (b\delta+1)K_1$ as a proper spanning subgraph. Then $\rho'>\rho(K_{n-(b+1)\delta-1}\cup K_{\delta}\cup (b\delta+1)K_1)=  n-(b+1)\delta-2$. Putting (\ref{equ::4}) into (\ref{equ::5}), and considering that $\rho'> n-(b+1)\delta-2$, we have
\begin{equation} \label{equ::6}
      y_3=\frac{\rho' y_2}{\rho'-(n-(b+1)\delta-2)}.
\end{equation}
Recall that $n\geq (b+1)s+2$. Then $\delta+1\leq s\leq(n-2)/(b+1)$. Since $G_{s}^2$ is not a regular graph, it follows that $\rho<n-1$. Suppose to the contrary that $\rho\geq\rho'$. Combining this with (\ref{equ::3}) and (\ref{equ::6}) yields that
\begin{equation*}
\begin{aligned}
   & y^{T}(\rho'-\rho)x\\
   &=y^{T}(A(G'_s)-A(G_{s}^2))x\\
   &=\sum_{i=b\delta\!+\!2}^{bs+1}\sum_{j=1}^{n\!-\!(b\!+\!1)s\!-\!1}\!(x_{u_{i}}y_{w_{j}}\!+\!x_{w_{j}}y_{u_{i}})\!
+\sum_{i=b\delta\!+\!2}^{bs}\sum_{j=i\!+\!1}^{bs\!+\!1}\!(x_{u_{i}}y_{u_{j}}\!+\!x_{u_{j}}y_{u_{i}})\! \!-\!
\sum_{i=\delta\!+\!1}^{s}\sum_{j=1}^{b\delta\!+\!1}\!(x_{v_{i}}y_{u_{j}}\!+\!x_{u_{j}}y_{v_{i}})\!\\
   &=(s\!-\!\delta)[b(n\!-\!(b\!+\!1)s\!-\!1)(x_{2}y_{3}\!+\!x_{3}y_{3})\!+\!b(bs\!-\!b\delta\!-\!1)x_{2}y_{3}\!-\!(b\delta\!+\!1)(x_{1}y_{2}\!+\!x_{2}y_{3})]\\
   &=(s\!-\!\delta)[(b(n\!-\!(b\!+\!1)\delta\!-\!s\!-\!2)\!-\!1)x_{2}y_{3}\!+\!b(n\!-\!(b\!+\!1)s\!-\!1)x_{3}y_{3}\!-\!(b\delta\!+\!1)x_{1}y_{2}]\\
   &\geq(s\!-\!\delta)[b(n\!-\!(b\!+\!1)\delta\!-\!s\!-\!2)x_{2}y_{3}\!-\!(b\delta\!+\!1)x_{1}y_{2}]  ~(\mbox{since $x_{3}\!\geq\! x_{2}$ and $n\geq (b\!+\!1)s\!+\!2$})\\
   &\geq(s\!-\!\delta)[b(n\!-\!(b\!+\!1)\delta\!-\!\frac{n\!-\!2}{b\!+\!1}\!-\!2)x_{2}y_{3}\!-\!(b\delta\!+\!1)x_{1}y_{2}] ~~(\mbox{since $s\leq\frac{n-2}{b+1}$})\\
   &=(s\!-\!\delta)x_{1}y_{2}\left(\frac{bs\rho'(n-(b+1)\delta-\frac{n-2}{b+1}-2)}{\rho(\rho'-(n-(b+1)\delta-2))}-(b\delta+1)\right)\\
   &\geq\frac{(b\delta\!+\!1)(s\!-\!\delta)x_{1}y_{2}}{\rho(\rho'\!-\!(n\!-\!(b\!+\!1)\delta\!-\!2))}[\rho'(n\!-\!(b\!+\!1)\delta\!-\!
                    \frac{n\!-\!2}{b\!+\!1}\!-\!2)\!-\!\rho(\rho'\!-\!(n\!-\!(b\!+\!1)\delta
                    \!-\!2))]\\
   &(\mbox{since $s\geq \delta+1$ and $b\geq 1$})\\
   &\geq\frac{(b\delta\!+\!1)(s\!-\!\delta)x_{1}y_{2}}{\rho(\rho'\!-\!(n\!-\!(b\!+\!1)\delta\!-\!2))}(\rho'(2n\!-\!2(b\!+\!1)\delta\!-\!\frac{n\!-\!2}{b\!+\!1}\!-\!4)\!-\!\rho\rho')~~(\mbox{since $\rho\geq\rho'$})\\
    &\geq\frac{(b\delta\!+\!1)(s\!-\!\delta)\rho'(\frac{3n}{2}\!-\!2(b\!+\!1)\delta\!-\!3\!-\!\rho)x_{1}y_{2}}{\rho(\rho'-(n-(b+1)\delta-2))}~~(\mbox{since $b\geq 1$})\\
    &\geq\frac{(b\delta\!+\!1)(s\!-\!\delta)\rho'(n\!-\!1\!-\!\rho)x_{1}y_{2}}{\rho(\rho'\!-\!(n\!-\!(b\!+\!1)\delta\!-\!2))} ~~(\mbox{since $n\geq 4(b\!+\!1)\delta\!+\!4$})\\
    &>0~~(\mbox{since $\rho< n\!-\!1$, $\rho'>n\!-\!(b\!+\!1)\delta\!-\!2$ and $s\geq \delta\!+\!1$}).
\end{aligned}
\end{equation*}
This implies that $\rho'>\rho$, which contradicts the assumption $\rho'\leq\rho$. It follows that $\rho'>\rho$. Combining this with (\ref{equ::1}) and (\ref{equ::2}), we may conclude that $\rho(G)\leq \rho(G_s^1)\leq \rho(G_s^2)<\rho(T(n,b,\delta))$.

{\flushleft\bf{Case 2.}} $s<\delta$.

Let $G_{s}^{3}=K_{s} \nabla (K_{n-s-(\delta+1-s)(t-1)}\cup (t-1)K_{\delta+1-s})$. Recall that $G$ is a spanning subgraph of $G_s^1=K_{s} \nabla (K_{n_1}\cup K_{n_2}\cup \cdots \cup K_{n_t})$, where $t=bs+2$, $n_1\geq n_2\geq \cdots\geq n_t$ and  $\sum_{i=1}^tn_i=n-s$. Clearly, $n_t\geq \delta+1-s$ because the minimum degree of $G_s^1$ is at least $\delta$. By Lemma \ref{lem::3.2}, we have
 \begin{equation}\label{equ::7}
\rho(G_{s}^{1})\leq \rho(G_{s}^{3}),
\end{equation}
where the equality holds if and only if $(n_1,\ldots,n_t)=(n-s-(\delta+1-s)(t-1),\delta+1-s,\ldots,\delta+1-s)$.

Assume that $\rho(G^{3}_{s})=\rho^*\geq n -(\delta+1-s)(t-1)$.
Let $x$ be the Perron vector
of $A(G^{3}_{s})$ corresponding to $\rho^*$. By symmetry, $x$ takes the same values  $x_{1}$, $x_{2}$, and $x_{3}$ on the vertices of $K_{n-s-(\delta+1-s)(t-1)}$, $K_{\delta+1-s}$ and $K_s$, respectively. Then, by $A(G^{3}_{s})x=\rho^* x$, we have
\begin{eqnarray}
\rho^* x_{1}&\!\!=\!\!&(n-1-(\delta+1-s)(t-1)-s)x_{1}+sx_{3},\label{equ::8}\\
\rho^*x_{2}&\!\!=\!\!&(\delta-s)x_{2}+sx_{3}, \label{equ::9}\\
\rho^*x_{3}&\!\!=\!\!&(n\!-\!(\delta\!+\!1\!-\!s)(t\!-\!1)\!-\!s)x_{1}\!+\!
    (\delta\!+\!1\!-\!s)(t\!-\!1)x_{2}\!+\!(s\!-\!1)x_{3}.\label{equ::10}
\end{eqnarray}
From (\ref{equ::8}) and (\ref{equ::9}), we obtain
\begin{equation}\label{equ::11}
\left\{
\begin{aligned}
 &x_{1}=\frac{sx_{3}}{\rho^*-(n-1-(\delta+1-s)(t-1))+s} , \\
 &x_{2}=\frac{sx_{3}}{\rho^*-\delta+s} .
   \end{aligned}
   \right.
\end{equation}
Note that $n\geq F(b,\delta)\geq b\delta^3+\delta$, $b\geq 1$ and $\delta\geq s+1\geq 2$. Then $\rho^*\geq n-(\delta+1-s)(t-1)> \delta+1$. Putting (\ref{equ::11}) into (\ref{equ::10}), and considering that $\rho^*\geq n-(\delta+1-s)(t-1)>\delta+1$, we get
\begin{equation*}
\begin{aligned}
                         &~~~~ \rho^*+1\\
                         &=s+\frac{s(n-(\delta+1-s)(t-1)-s)}{\rho^*-(n-1-(\delta+1-s)(t-1))+s}+
                          \frac{s(\delta+1-s)(t-1)}{\rho^*-(\delta-s)}\\
                         &< s+\frac{s(n-(\delta+1-s)(t-1)-s)}{s+1}+\frac{s(\delta+1-s)(t-1)}{s+1}\\
                         &=s+\frac{s(n-s)}{s+1}\\
                          &=n\!-\!(\delta\!+\!1\!-\!s)(\!t-\!1)\!-\!\frac{n\!-\!s\!-\!(\delta\!+\!1\!-\!s)(t\!-\!1)(s\!+\!1)}{s\!+\!1}\\
                               &\leq n\!-\!(\delta\!+\!1\!-\!s)(t\!-\!1)\!-\!\frac{b\delta^{3}\!+\!\delta\!-\!s\!-\!(\delta
                         \!+\!1\!-\!s)(t\!-\!1)(s\!+\!1)}{s\!+\!1} ~~(\mbox{since ~$n\geq b\delta^{3}\!+\!\delta$})\\
                               &<n\!-\!(\delta\!+\!1\!-\!s)(t\!-\!1)~~(\mbox{since ~$\delta\geq s+1$, $b\geq 1$ and $s\geq 1$})\\
                                 &\leq\rho^*,
 \end{aligned}
\end{equation*}
which is impossible. Thus we have
\begin{equation}\label{equ::12}
\begin{aligned}
                      \rho^*&<n\!-\!(\delta\!+\!1\!-\!s)(t\!-\!1)\\
                         &=n\!-\!b\delta\!-\!2\!-\!((s\!-\!1)(\delta\!-\!s)b\!+\!\delta\!-\!(s\!+\!1))~~(\mbox{since ~$\delta\geq s\!+\!1$, $b\!\geq\!1$ and $s\!\geq\!1$})\\
                         &\leq n\!-\!b\delta\!-\!2.
\end{aligned}
\end{equation}
Since $T(n,b,\delta)$ contains $K_{n-b\delta-1}\cup (b\delta+1)K_{1}$ as a proper spanning  subgraph, it follows that $\rho(T(n,b,\delta))>\rho(K_{n-b\delta-1}\cup (b\delta+1)K_{1})=n-b\delta-2$. Combining this with (\ref{equ::1}), (\ref{equ::7}) and (\ref{equ::12}), we may conclude that $\rho(G)\leq \rho(G_s^1)\leq \rho(G_s^3)<\rho(T(n,b,\delta))$.

{\flushleft\bf{Case 3.}} $s=\delta.$

Note that $G_{s}^1$ is a spanning subgraph of $T(n,b,\delta)$. Then
$$\rho(G_{s}^1)\leq \rho(T(n,b,\delta)),$$
with equality holding if and only if $G_{s}^1\cong T(n,b,\delta)$. Combining this with (\ref{equ::1}), we may conclude that
$$\rho(G)\leq\rho(T(n,b,\delta)),$$
where the equality holds if and only if $G\cong T(n,b,\delta)$. This completes the proof. \end{proof}

From Theorem \ref{thm::1.2}, we can get the result of Liu, Liu and Feng\cite{W.L} immediately.
\begin{cor}(See \cite{W.L})\label{cor::3.1}
Let $G$ be a connected graph of even order $n$ with minimum degree $\delta(G)\geq 2$. If $n\geq \max\{7+7\delta+2\delta^{2}, \delta^{3}+3\delta^{2}+2\delta\}$ and $\rho(G)\geq\rho(T(n,1,\delta))$,
then $G$ has a perfect matching, unless $G\cong T(n,1,\delta)$.
\end{cor}

\section{Proof of Theorem \ref{thm::1.3}}
The following two structural lemmas will play an essential role in the proof of Theorem \ref{thm::1.3}.
\begin{lem}(See \cite{Y.L})\label{lem::4.1}
Let $G$ be a graph of order $n$ with minimum degree $\delta(G)$ and let $a$ and $b$ be integers such that $1\leq a< b$. Then if $\delta(G)\geq a$, $n\geq 2a+b+\frac{a^{2}-a}{b}$ and
$$\max\{d_{G}(u), d_{G}(w)\}\geq\frac{an}{a+b},$$
for any two nonadjacent vertices $u$ and $w$ of $G$, then $G$ contains an $[a,b]$-factor.
\end{lem}

\begin{lem}(See \cite{T.N})\label{lem::4.2}
Suppose that $k\geq 3$. Let $G$ be a connected graph of order $n\geq 4k-3$ with minimum degree $\delta(G)$, where $k\cdot n$ is even and $\delta(G)\geq k$. If
$$\max\{d_{G}(u), d_{G}(w)\}\geq\frac{n}{2},$$
for any two nonadjacent vertices $u$ and $w$ of $G$, then $G$ contains a $k$-factor.
\end{lem}
\begin{lem}(See \cite{M.F})\label{lem::4.3}
Let $G$ be a connected graph of order $n$. If
\begin{equation}\label{equ::13}
\begin{aligned}
\rho(G)\geq n-2,
\end{aligned}
\end{equation}
then $G$ contains a Hamiltonian path unless $G\cong H_{n,1}$. If the inequality (\ref{equ::13}) is
strict, then $G$ contains a Hamiltonian cycle unless $G\cong H_{n,2}$.
\end{lem}

\begin{lem}(See \cite{BH})\label{lem::4.4}
Let $M$ be a real symmetric matrix, and let $\lambda_{1}(M)$ be the largest eigenvalue of $M$. If $B_\Pi$ is an equitable quotient matrix of $M$, then the eigenvalues of  $B_\Pi$ are also eigenvalues of $M$. Furthermore, if $M$ is nonnegative and irreducible, then $\lambda_{1}(M) = \lambda_{1}(B_\Pi).$
\end{lem}

\begin{lem}\label{lem::4.5}
Suppose that $G$ is a connected graph of order $n$. Let $u,w$ be two nonadjacent vertices of $G$ such that $\max\{d_{G}(u), d_{G}(w)\}\leq t$, where $t\geq 1$. Then $\rho(G)\leq \rho(K_{t} \nabla (2K_{1}\cup K_{n-t-2}))$, with equality if and only if $G\cong K_{t} \nabla (2K_{1}\cup K_{n-t-2})$.
\end{lem}
\renewcommand\proofname{\bf Proof}
\begin{proof}
Let $x$ be the Perron vector of
$A(G)$ with respect to $\rho(G)$. Assume that $V(G)\setminus\{u,w\}=\{v_{1},v_{2},\ldots,v_{n-2}\}$ with $x_{v_1}\geq x_{v_2}\geq\cdots\geq x_{v_{n-2}}$. Let
$$G'=G-\sum_{v\in N_{G}(u)}uv-\sum_{v\in N_{G}(w)}wv+\sum_{i=1}^{t}(uv_{i}+wv_{i}).$$
Since $d_{G}(u)\leq t$ and $d_{G}(w)\leq t$, by Lemma \ref{lem::2.1}, it follows that $\rho(G)\leq \rho(G')$, where equality holds if and only if $G\cong G'$. Obviously, $G'$ is a spanning subgraph of $K_{t} \nabla (2K_{1}\cup K_{n-t-2})$. Thus
\begin{equation*}
\begin{aligned}
\rho(G)\leq \rho(K_{t} \nabla (2K_{1}\cup K_{n-t-2})),
\end{aligned}
\end{equation*}
with equality if and only if $G\cong K_{t} \nabla (2K_{1}\cup K_{n-t-2})$.
\end{proof}

Now we shall give the proof of Theorem \ref{thm::1.3}.
\renewcommand\proofname{\bf Proof of Theorem \ref{thm::1.3}}
\begin{proof}
Let $G$ be a graph of order $n\geq 3a+b-1$ with $\rho(G)>\rho(H_{n,a})$, where $a\cdot n$ is an even integer and $b\geq a\geq1$. We first assert that $G$ is connected. If not, suppose that $G_{1}, \ldots, G_{p}$ ($p\geq 2$) are the components of $G$. Then $\rho(G)=\max\{\rho(G_{1}),\ldots,\rho(G_{p})\}\leq \rho (K_{n-1})=n-2$, which contradicts that $\rho(G)>\rho(H_{n,a})\geq n-2$, and so we assume that $G$ is connected in the following. We next assert that $\delta(G)\geq a$. Otherwise, let $1 \leq\delta(G)\leq a-1$. Then $G$ is a spanning subgraph of $H_{n,a}$ for $a\geq 2$, and so $\rho(G)\leq \rho(H_{n,a})$, a contradiction.

If $a=b=1$. Since $n$ is even and $\rho(G)>\rho(H_{n,1})$, by Lemma \ref{lem::4.3}, it follows that $G$ contains a Hamiltonian path of even number, and so $G$ contains a 1-factor. If $a=b=2$, then by using the same analysis as above, we deduce that $G$ contains a 2-factor. Next, we assume that both $a$ and $b$ are not equal to 1 or 2 simultaneously. Suppose that $G$ contains no $[a,b]$-factor of order $n\geq 3a+b-1$. By Lemmas \ref{lem::4.1} and \ref{lem::4.2}, there exist two nonadjacent vertices $u$ and $w$ such that $\max\{d_{G}(u), d_{G}(w)\}\leq \frac{an}{a+b}-1$.
Let $t=\frac{an}{a+b}-1$. Then $t\geq \delta(G)\geq a\geq1$. Suppose that $G_{1}\cong K_{t} \nabla (2K_{1}\cup K_{n-t-2})$. By Lemma \ref{lem::4.5}, we obtain
\begin{equation}\label{equ::14}
\begin{aligned}
\rho(G)\leq \rho(G_{1}).
\end{aligned}
\end{equation}
The vertex set of $G_{1}$ can be partitioned as $V(G_1)=V(2K_1)\cup V(K_t)\cup V(K_{n-t-2})$, where $V(2K_1)=\{u, w\}$, $V(K_{t})=\{v_{1},\ldots,v_{t}\}$ and $V(K_{n-t-2})=\{v_{t+1},\ldots,v_{n-2}\}$. Let $G'_1=G_{1}-\{uv_{a},uv_{a+1},\ldots,uv_{t}\}+\{wv_{t+1},wv_{t+2},\ldots,wv_{n-2}\}$. Clearly,
$G'_1\cong H_{n,a}$. If $a=1$ and $b\geq 2$, then $G'_1\cong K_{1}\cup K_{n-1}$. Note that $t=\frac{an}{a+b}-1$. Then $n\geq 3t+3$. Observe that $A(G_1)$ has the equitable quotient matrix
$$
B_1=\begin{bmatrix}
0 &t &  0\\
2&t-1 & n-t-2\\
  0&t& n-t-3\\
\end{bmatrix}.
$$
The characteristic polynomial of $B_1$ is
\begin{equation*}
      f(\lambda)=\lambda^3-(n-4)\lambda^2-(n+2t-3)\lambda+2tn-2t^2-6t.
\end{equation*}
By a simple computation, we have
\begin{equation*}
\begin{aligned}
 f(-\infty)=-\infty,~~ f(0)=2t(n-t-3)>0,~~ f(n-3)=-2t^2<0,
\end{aligned}
\end{equation*}
and
\begin{equation*}
\begin{aligned}
 f(n-2)=(n-\frac{3}{2})^2-2t^2-2t-\frac{1}{4}\geq 7t^2+7t+2>0.
\end{aligned}
\end{equation*}
By Lemma \ref{lem::4.4}, we have $\lambda_{1}(B_1)=\rho(G_1)$. Notice that $\rho(G'_1)=\rho(K_{1}\cup K_{n-1})=n-2$ and $\lambda_{1}(B_1)<n-2$. Then $\rho(G_1)<\rho(G'_1)$. Combining this with (\ref{equ::14}), we obtain that $\rho(G)\leq\rho(G_1)<\rho(H_{n,1})$, a contradiction. Next, we assume that $a\geq 2$ and $b\geq3$. Recall that $t=\frac{an}{a+b}-1$. Then $n\geq 2t+2$. Let $x$ be the Perron vector of $A(G_1)$ with respect to $\rho(G_1)$, and let $\rho=\rho(G_1)$. By symmetry, $x$ takes the same value (say $x_1$, $x_2$ and $x_3$) on the vertices of $V(2K_1)$, $V(K_t)$ and $V(K_{n-t-2})$, respectively. Then, by $A(G_1)x=\rho x$, we have
 \begin{equation}\label{equ::15}
x_{2}=\frac{(\rho-(n-t-3))x_3}{t}.
\end{equation}
Similarly, let $y$ be the Perron vector of $A(G'_1)$ corresponding to $\rho(G'_1)=\rho'$. By symmetry, $y$ takes the same values $y_{1}$, $y_{2}$ and $y_{3}$ on the vertices of $V(K_1)$, $V(K_{a-1})$ and $V(K_{n-a})$, respectively. Then, by $A(G'_1)y=\rho' y$, we obtain
\begin{eqnarray}
 &\rho' y_1&=(a-1) y_{2},\label{equ::16}\\
 &\rho' y_3 &=(a-1) y_{2}+(n-a-1) y_{3}.\label{equ::17}
\end{eqnarray}
Since $G'_{1}$ contains $K_{1}\cup K_{n-1}$ as a proper spanning subgraph, we have $\rho'> n-2 > n-a-1$. Putting (\ref{equ::16}) into (\ref{equ::17}), and considering that $\rho'> n-a-1$, we have
\begin{equation} \label{equ::18}
      y_3=\frac{\rho' y_1}{\rho'-(n-a-1)}.
\end{equation}
Obviously, neither $G_1$ nor $G'_1$ is a regular graph. Thus $\rho<n-1$ and $\rho'<n-1$. Combining this with (\ref{equ::15}) and (\ref{equ::18}) yields that
\begin{equation*}
\begin{aligned}
  y^{T}(\rho'-\rho)x&= y^{T}(A(G'_1)-A(G_1))x\\
                    &=\sum_{i=t+1}^{n-2}(x_{w}y_{v_{i}}+x_{v_{i}}y_{w})-\sum_{i=a}^{t}(x_{u}y_{v_{i}}+x_{v_{i}}y_{u})\\
                    &=(n\!-\!t\!-\!2)(x_{1}y_{3}+x_{3}y_{3})\!-\!(t\!-\!a+1)(x_{1}y_{3}+x_{2}y_{1})\\
                    &=(n-2t+a-3)x_{1}y_{3}+(n-t-2)x_{3}y_{3}-(t-a+1)x_{2}y_{1}\\
                    &>(n-t-2)x_{3}y_{3}-(t-a+1)x_{2}y_{1} ~~(\mbox{since $n\geq 2t+2$ and $a\geq 2$}) \\
                    &=x_{3}y_{1}\left(\frac{\rho'(n-t-2)}{\rho'-(n-a-1)}-\frac{(\rho-(n-t-3))(t-a+1)}{t}\right)\\
                     &>x_{3}y_{1}\left(\frac{\rho'(n\!-\!t\!-\!2)}{a}\!-\!\frac{(t\!-\!a\!+\!1)(t+2)}{t}\right)~~(\mbox{since $\rho\!<\!n\!-\!1$ and $\rho'\!<\!n\!-\!1$})\\
\end{aligned}
\end{equation*}
\begin{equation*}
\begin{aligned}
                    &> x_{3}y_{1}\left(\frac{(n-2)(n-t-2)}{a}-\frac{(t-a+1)(t+2)}{t}\right)~~(\mbox{since $\rho'>n-2$})\\
                    &\geq \frac{((n-2)(n-t-2)-(t-a+1)(t+2))x_{3}y_{1}}{t}~(\mbox{since $t\geq a$})\\
                    &>0 ~~(\mbox{since $t\geq a\geq2$ and $n\geq 2t+2$}).
\end{aligned}
\end{equation*}
Therefore, we have $\rho< \rho'$. Combining this with (\ref{equ::14}), we may conclude that $\rho(G)\leq\rho(G_1)<\rho(H_{n,a})$, a contradiction. Therefore, the proof is complete.
\end{proof}

\section{Concluding remarks}
A \textit{fractional matching} of a graph $G$ is a function $f$ giving each edge a number in $[0,1]$ such that $\sum_{e\in \Gamma(v)}f(e)\leq 1$ for each $v\in V(G)$, where $\Gamma(v)$ is the set of edges incident to $v$. A \textit{fractional perfect matching} of a graph $G$ is a fractional matching $f$ with $\sum_{e\in E(G)}f(e)=\frac{n}{2}$. The relations between the eigenvalues and the fractional matching number of a graph were studied by several researchers \cite{O2,PLZ,XZS}.

The following fundamental lemmas provide some sufficient and necessary conditions for the existence of a $[1,b]$-factor for $b\geq 2$ and a fractional perfect matching in a graph, respectively. Let $i(H)$ be the number of isolated vertices of a graph $H$.
\begin{lem}(See \cite{C.B})\label{lem::5.1}
Let $G$ be a graph and let $b\geq 2$ be a positive integer. Then $G$ contains a $[1,b]$-factor if and only if for every subset $S\subseteq V(G)$,
$$i(G-S)\leq b|S|.$$
\end{lem}

\begin{lem}(See \cite{E.S})\label{lem::5.2}
A graph $G$ has a fractional perfect matching if and only if for every subset $S\subseteq V(G)$,
$$i(G-S)\leq |S|.$$
\end{lem}
By Lemmas \ref{lem::5.1} and \ref{lem::5.2} and using the same analysis as the proof of Theorem \ref{thm::1.2}. We easily obtain the following result.
\begin{thm}\label{thm::5.1}
Suppose that $G$ is a connected graph of order $n$ with minimum degree $\delta$.
\begin{enumerate}[(i)]
\item If $n\geq 4(b+1)\delta+4$ and $\rho(G)\geq\rho(K_{\delta} \nabla (K_{n-(b+1)\delta-1}\cup (b\delta+1)K_{1}))$ where $b\geq 2$, then $G$ contains a $[1,b]$-factor, unless $G\cong K_{\delta} \nabla (K_{n-(b+1)\delta-1}\cup (b\delta+1)K_{1})$.
\item If $n\geq 8\delta+4$ and $\rho(G)\geq\rho(K_{\delta} \nabla (K_{n-2\delta-1}\cup (\delta+1)K_{1}))$, then $G$ contains a fractional perfect matching, unless $G\cong K_{\delta} \nabla (K_{n-2\delta-1}\cup (\delta+1)K_{1})$.
\end{enumerate}
\end{thm}
Corollary 1.1 provides a spectral condition to guarantee the existence a $k$-factor
in a graph. It is natural to ask the following question for $k\geq 1$.

\begin{prob}\label{prob::5.1}
What is the maximum spectral radius and what is the corresponding extremal graph among all graphs with a unique $k$-factor for $k\geq 1$?
\end{prob}

In this paper, we give the answer to Problem \ref{prob::5.1} for $k=1$. However, the structure of graphs with a unique $k$-factor is more complicated for $k\geq 2$, and it seems difficult to determine the extremal graphs about the problem.

If $k\leq n$. Let $n=sk+t$ with $s\geq 1$ and $0\leq t\leq k-1$. We give the process to construct the graph $G(2n,k)$. First define a graph $F_1$ on $2(k+t)$ vertices as follows. Let $H_{1}\cong K_{t}\nabla tK_1$, and let $A_{11}=V(tK_{1})$ and $A_{12}=V(K_t)$. Denote by $H_2$ the graph obtained from $kK_{1}\cup K_{k}$ by adding edges between $V(kK_{1})$ and $V(K_{k})$ such that $d_{H_2}(v)=k-t$ for $v\in V(kK_{1})$ and $d_{H_2}(u)=2k-t-1$ for $u\in V(K_{k})$. Suppose that $A_{21}=V(kK_{1})$ and $A_{22}=V(K_{k})$. Let $F_1$ be the graph of order $2(k+t)$ obtained from $H_{1}\cup H_{2}$ by connecting all vertices of $A_{1j}$ with $A_{2(3-j)}$ for $1\leq j\leq 2$ and adding all edges between $A_{12}$ and $A_{22}$. The resulting graph $F_1$ has exactly one $k$-factor. Suppose that $U_1=A_{11}\cup A_{21}$ and $W_1=A_{12}\cup A_{22}$. Next take $s-1$ copies of $K_{k}\nabla kK_1$ labeled $F_{2},\ldots, F_{s}$. For $2\leq i\leq s$, let $U_{i}$ and $W_{i}$ be the vertices set of $V(kK_1)$ and $V(K_k)$ in each $F_i$, respectively. Then the graph $G(2n,k)$ is obtained by adding edges connecting all vertices of $W_i$ in $F_i$ to all vertices in $F_j$ for each $i,j$ with $1\leq i< j\leq s$. The resulting graph $G(2n,k)$ has a unique $k$-factor. We end the paper by proposing the following problem for further research.

\begin{prob}
For $k\geq 1$. Suppose that $G$ is a graph of order $2n$ with a unique $k$-factor.

\noindent (I) Does $\rho(G)\leq\rho(K_{2n-k}\nabla H)$, where $H$ is a $2(k-n)$-regular graph if $k> n$?

\noindent (II) Does $\rho(G)\leq\rho (G(2n,k))$, if $k\leq n$?
\end{prob}

\end{document}